\pdfoutput=1
\documentclass[10pt,a4paper]{article}

\usepackage{arxiv}
\usepackage[utf8]{inputenc}
\usepackage[T1]{fontenc}
\usepackage{amsmath,amssymb,amsfonts}
\usepackage{graphicx}
\usepackage{booktabs}
\usepackage{tikz}
\usetikzlibrary{arrows.meta,positioning,decorations.pathreplacing,calc,shapes,fit,backgrounds}
\usepackage{float}
\usepackage[numbers,sort&compress]{natbib}
\usepackage{hyperref}
\usepackage[nameinlink,capitalize]{cleveref}

\graphicspath{{figures/}}

\newcommand{\Dc}{{}^{C}\!D}
\newcommand{\Drl}{{}^{RL}\!D}

\title{Fixed-grid schemes for fractional moving-boundary problems}

\author{%
  Conrard Giresse Tetsassi Feugmo\textsuperscript{1,2,*}\\[4pt]
  \small\textsuperscript{1}Department of Chemistry, University of Waterloo,
  Waterloo, Ontario N2L 3G1, Canada\\
  \small\textsuperscript{2}Department of Physics and Astronomy, University of Waterloo,
  Waterloo, Ontario N2L 3G1, Canada\\[2pt]
  \small\textsuperscript{*}Corresponding author: \href{mailto:cgtetsas@uwaterloo.ca}{cgtetsas@uwaterloo.ca}
}

\date{}

\begin{document}

\maketitle

\begin{abstract}
The Landau transformation $\xi=x/s(t)$ is the usual way to hold a moving
boundary in place. It fails when the equation has a fractional time
derivative, because fractional derivatives have no chain rule: the transformed
operator gains an extra nonlocal integral that has a closed form only when the
front is self-similar. We draw two consequences. First, we build a
sharp-interface finite-difference scheme for the conservation-consistent
(memory-flux) fractional Stefan problem on a fixed grid, storing the memory
limit $h(x)=s^{-1}(x)$ as the time each node melted. The front is advanced by
an integrated energy balance, because the local gradient at the front is zero,
and the singular source is timed to the sub-step crossing instant. Against a
known self-similar solution the front constant is accurate to
$0.02$--$0.11\,\%$, while the competing formulations differ by
$17$--$27\,\%$ at $\alpha=0.6$. A remaining positivity violation follows from
the scaling $\Delta t/t_0\sim s_0^{1-2/\alpha}$ and disappears on a
front-locked time grid. We also derive two exact solutions whose fronts are
not powers of $t$, one in the source-dominated limit and a Mittag-Leffler
traveling wave for the full problem, and verify the scheme against them.
Second, neural-network solvers for moving boundaries are often trained in a
normalised coordinate; for $\alpha<1$ this silently changes the problem, and
the change vanishes at $\alpha=1$, where checks are usually run. In the
classical limit our scheme is more accurate in $12$\,ms than a reimplemented
moving-boundary network in $177$\,s.
\end{abstract}

\textbf{Keywords:} fractional Stefan problem, moving boundary, memory flux, front fixing, physics-informed neural networks, verification

\section{Introduction}

Moving-boundary problems are usually attacked in one of two ways~\citep{Crank1984}: capture the interface on a fixed grid, or fix the interface
by transforming the domain. The second is the Landau transformation
$\xi=x/s(t)$~\citep{Landau1950}, which maps
the growing domain to $[0,1]$ at the cost of an advective term, and it is the
default in one dimension because it makes the free boundary an ordinary unknown.

This paper starts from one observation: the Landau transformation does not
work when the time derivative is fractional. Section~\ref{sec:nocr}
states the obstruction. Sections~\ref{sec:form}--\ref{sec:noneq} build and verify
a fixed-grid scheme that does not need it. Section~\ref{sec:pinn} turns to the
classical limit, where the transformation \emph{is} legitimate, uses the
resulting solver to benchmark a recent moving-boundary neural network, and
explains why the same obstruction applies to neural solvers.
Sections~\ref{sec:limit}--\ref{sec:fix} find the cause of the scheme's one
remaining defect, a small positivity violation, and remove it.

Verifying the scheme also forces a solution worth stating as a result in its
own right. To our knowledge, every exact solution to the fractional Stefan
problem with the standard Caputo derivative published to date, including those collected in the survey of~\citet{RogosinDubatovskaya2023}, has a free boundary that moves as a power of
time,
\[
s(t)\;\propto\;t^{\alpha/2},
\]
so that a single similarity variable collapses the governing equation to an
ordinary one. Section~\ref{sec:tw} constructs a solution on the opposite
mechanism instead: a front of \emph{constant} speed,
\[
s(t)=s_0+vt,
\]
whose profile depends on $x$ and $t$ only through the time elapsed since
melting at that point. A front that translates rather than spreads admits no
similarity variable, and the closed form is correspondingly not a Wright
function but the Mittag-Leffler function
\[
W(\tau)=E_{2-\alpha}\bigl(v^{2}\tau^{\,2-\alpha}\bigr).
\]
This is what lets the scheme be verified against a front that no self-similar
benchmark can test.

\section{The obstruction}\label{sec:nocr}

\begin{figure}[t]
\centering
\includegraphics[width=\linewidth]{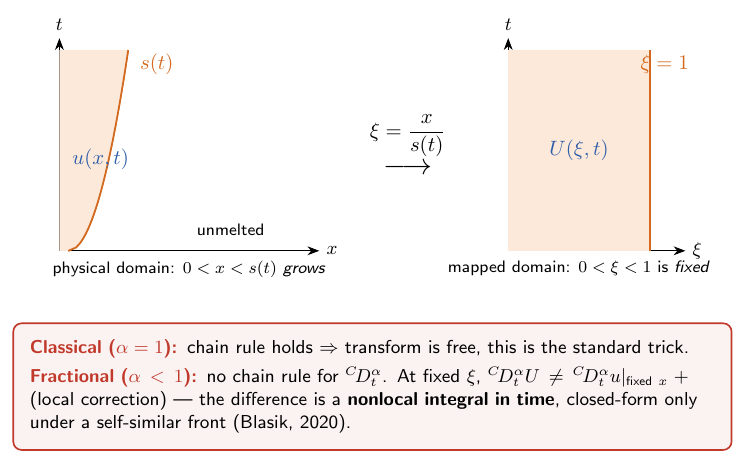}
\caption{The Landau transformation $\xi=x/s(t)$ fixes the growing domain
$0<x<s(t)$ (left) onto $0<\xi<1$ (right) at the cost of an advective term.
Classically ($\alpha=1$) this is free: the chain rule writes the transformed
time derivative as the original plus a local correction. Fractionally
($\alpha<1$) there is no chain rule for $\Dc^\alpha_t$: the difference between
the derivative at fixed $\xi$ and at fixed $x$ is a nonlocal integral in time,
reducible to closed form only when the front is self-similar~\citep{Blasik2020}.}
\label{fig:landau}
\end{figure}

There is no chain rule for fractional derivatives~\citep{OldhamSpanier1974,Tarasov2016} (Fig.~\ref{fig:landau}). A Caputo
operator taken at fixed $\xi=x/s(t)$ is therefore not the operator at fixed $x$
plus pointwise corrections: the difference is an additional nonlocal integral in
time. \citet{Blasik2020} writes this correction down explicitly, and can resum it
only because the front has already been assumed to be $s\propto\tau^{\alpha/2}$,
which turns the transformation Jacobian into a power of $\tau$. For a general
front (the case whenever the boundary datum varies, the medium is finite, or the
initial front is not a point), no such reduction exists.

Front-fixing schemes for fractional moving-boundary problems~\citep{BlasikKlimek2015,Blasik2020,Blasik2020b} are consequently restricted to reproducing a
similarity solution that has been put in by hand. Anything else requires a fixed
grid with the domain history represented explicitly, which is what
Section~\ref{sec:scheme} constructs.

\section{Formulations}\label{sec:form}

Write the one-phase problem as a balance law with a constitutive law,
\begin{equation}
\partial_t H = -\partial_x q,\qquad q=-k\,\partial_x u,
\label{eq:balance}
\end{equation}
with $H=u+\mathbf{1}_{x<s(t)}$ the enthalpy. The common fractional generalisation
substitutes a Caputo derivative into the bulk equation \emph{and} into the
classical Stefan condition,
\begin{equation}
\Dc^\alpha u = u_{xx},\qquad \Dc^\alpha s(t)=-\lambda\,u_x(s(t),t),
\label{eq:naive}
\end{equation}
which admits Wright-function similarity solutions~\citep{RoscaniSantillan2014};
their front prefactor solves a transcendental equation, for which
\citet{Caruso2025} give efficient numerical schemes. \citet{KubicaRyszewska2021} note that these
conditions ``are not the consequences of the assumptions imposed on the flux''.
Leaving \eqref{eq:balance} intact and putting the memory in the constitutive law,
\begin{equation}
q(x,t) = -\Drl^{\,1-\alpha}_{h(x)}\,\partial_x u(x,t),
\qquad h(x)=s^{-1}(x),
\label{eq:memflux}
\end{equation}
gives instead~\citep{KubicaRyszewska2021,RoscaniBollatiTarzia2018}
\begin{equation}
\Dc^\alpha_{h(x)}u-u_{xx}=
\begin{cases}
0, & x<s(0),\\[2pt]
-\dfrac{(t-h(x))^{-\alpha}}{\Gamma(1-\alpha)}, & s(0)<x<s(t),
\end{cases}
\label{eq:consistent}
\end{equation}
with $u(0,t)=\gamma$, $u(s(t),t)=0$, the interface condition
\begin{equation}
\dot s(t)=-\frac{1}{\Gamma(\alpha)}\lim_{y\to s(t)^-}\frac{\mathrm{d}}{\mathrm{d}t}
\int_{h(y)}^{t}(t-r)^{\alpha-1}u_x(y,r)\,\mathrm{d}r,
\label{eq:stefan}
\end{equation}
and an \emph{emergent} extra condition $u_x^-(s(t),t)=0$.

Two features of \eqref{eq:consistent} are absent from \eqref{eq:naive} and are
what make it conservative: the memory limit is $h(x)$, not $0$, so a material
point remembers only from the instant it melted (the only defensible reading on
a growing domain, since a fixed limit integrates over times at which the point
lay outside it), and the latent heat appears as a distributed source rather than
being carried entirely by the interface~\citep{VollerFalciniGarra2013}.

\citet{Ceretani2020} shows the competing formulations ``are nonequivalent though
all of them reduce to a classical Stefan problem when the order of the fractional
derivatives is replaced by one''. It follows that the $\alpha\to1$ limit cannot verify a scheme, because every
formulation passes it.

\subsection*{Relation to previous work}

\citet{VollerRoscani2023} derive \eqref{eq:consistent} from a thermodynamic
balance, note the zero-flux interface condition, and present an implicit
fixed-grid scheme using Gr\"unwald--Letnikov weights on a memory enthalpy;
\citet{RoscaniVoller2024} prove the memory-enthalpy and memory-flux formulations
equivalent and establish that the flux is continuous across the interface, the
formulation supplying a natural regularisation absent from the classical problem.
Together that work settles the modelling question and is the direct predecessor of
ours. The scheme of~\citet{VollerRoscani2023} regularises the interface through a
liquid fraction, so $h(x)$ enters only implicitly and the lower limit on the
enthalpy remains $0$; the front is captured rather than tracked; and the reported
computation is a single self-similar problem without a refinement study or a
comparison against an exact solution. \citet{RoscaniVoller2024} is analytical
throughout and reports no computation. We contribute a sharp-interface scheme with
$h(x)$ explicit and the emergent condition imposed, a source treatment correct to
sub-step resolution, and a verification. \citet{RogosinDubatovskaya2023} survey
the wider field. Other numerical approaches to Stefan-type problems include a
variational-iteration solution of the classical one-phase problem~\citep{Slota2007} and a finite-difference/spectral scheme for the
time-fractional diffusion equation~\citep{LinXu2007}, neither addressing the
front-fixing obstruction directly. \citet{ChmielowskaSlota2022} solve the
two-phase time-fractional Stefan problem on a fixed grid with the alternating
phase truncation method, an enthalpy approach in which the front is not
tracked explicitly; \citet{RoscaniTarzia2017} show that two
distinct fractional formulations both converge to the same classical limit, a
companion result to the nonequivalence away from that limit quantified in
Section~\ref{sec:noneq}.

\section{The scheme}\label{sec:scheme}

\subsection{Grid and activation times}

With $x_i=i\Delta x$ and $t_n=n\Delta t$, the front is captured between nodes and
we record for each node the exact instant $t^*_i$ at which it was crossed, by
linear interpolation of $s$ across the crossing step. Because $u(x_i,t)=0$ for
$t\le t^*_i$, increments before activation vanish identically and an L1 quadrature~\citep{Stynes2023}
with a fixed lower limit reproduces the variable-limit operator automatically.

\subsection{The front is advanced by an integrated balance}

\begin{figure}[t]
\centering
\includegraphics[width=\linewidth]{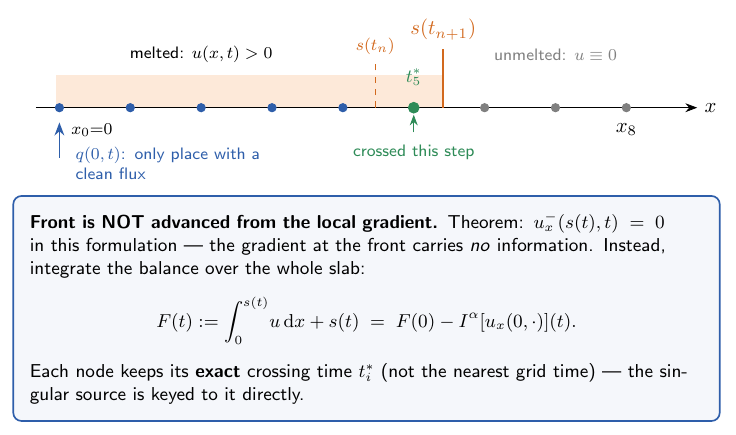}
\caption{The fixed spatial grid at a step $t_{n+1}$. Melted nodes ($x<s(t_{n+1})$,
shaded) carry $u>0$; unmelted nodes are exactly zero. Node $5$ was crossed
during this step and records its exact crossing time $t^*_5$ by linear
interpolation, rather than being rounded to the grid time $t_{n+1}$. The front
is not advanced from the local gradient, which \S\ref{sec:scheme} shows carries
no information about the front speed: it is advanced instead from the flux
$q(0,t)$ at the one fixed boundary where it is clean, via the integrated
balance~\eqref{eq:front}.}
\label{fig:gridscheme}
\end{figure}

That the local flux vanishes at the free boundary is a theorem, not an
artefact (Fig.~\ref{fig:gridscheme}):
if $\Dc^\alpha_{h(x)}u\in L^1(0,s(t))$ then $u_x^-(s(t),t)=0$~\citep{KubicaRyszewska2021,RoscaniVoller2024}, so the memory replaces the
classical flux jump by a continuous flux. The numerical consequence, which to our
knowledge has not been drawn, is that the local gradient carries no information
about the front speed and \eqref{eq:stefan} cannot be discretised directly. Integrating \eqref{eq:balance} over the slab, using $q=0$ ahead of the
front and $I^1\Drl^{1-\alpha}=I^\alpha$,
\begin{equation}
F(t):=\int_0^{s(t)}u\,\mathrm{d}x+s(t)=F(0)-I^{\alpha}\!\left[u_x(0,\cdot)\right](t).
\label{eq:front}
\end{equation}
The right-hand side involves only the flux at the \emph{fixed} boundary $x=0$,
where a one-sided second-order stencil is clean, and the moving boundary enters
through the smooth functional $\int_0^{s}u$. At $\alpha=1$ this is the classical
balance, so the classical and fractional solvers share one front update, a fact
Section~\ref{sec:pinn} uses.

\subsection{The singular source, and why sub-step timing matters}
\label{sec:source}

The source is singular at activation and must be integrated, not sampled:
\begin{equation}
S_i^{n}=\frac{(t_{n}-t^*_i)_+^{1-\alpha}-(t_{n-1}-t^*_i)_+^{1-\alpha}}
{\Delta t\,\Gamma(2-\alpha)} .
\label{eq:source}
\end{equation}
Replacing $t^*_i$ by the grid time collapses \eqref{eq:source} to
$\Delta t^{-\alpha}b_{m-1}/\Gamma(2-\alpha)$ with $b_j=(j+1)^{1-\alpha}-
j^{1-\alpha}$, the natural discretisation. It fails: the whole first-step impulse $b_0=1$ lands in one
cell, and since the increment a node receives as the front passes is
$\Delta u\approx-[(t_{n+1}-t^*_i)/\Delta t]^{1-\alpha}$, which tends to $-1$ as
$\alpha\to1$ for \emph{any} partial step, the field is driven negative, violating
the positivity established by~\citet{KubicaRyszewska2021}. The sub-step $t^*_i$
splits the impulse across two steps, reducing the violation three- to five-fold
while also improving the front constant (Table~\ref{tab:undershoot}).

\begin{table}[t]
\caption{Effect of sub-step source timing, $n_x=201$, $n_t=600$.
$\min u$ should be non-negative.}
\label{tab:undershoot}
\begin{tabular*}{\linewidth}{@{}llll@{}}
\toprule
$\alpha$ & $\min u$, integer index & $\min u$, sub-step & front error \\
\midrule
$0.70$ & $-0.078$ & $-0.015$ & $0.14\to0.02\,\%$ \\
$0.90$ & $-0.106$ & $-0.032$ & $0.33\to0.11\,\%$ \\
$0.95$ & $-0.119$ & $-0.035$ & --- \\
\bottomrule
\end{tabular*}
\end{table}

A residual violation of order $10^{-2}$ survives on the crossing step, at about
$35$ cells in $1.2\times10^5$. It is not a defect of the discretisation but the
finite-grid signature of a genuine singularity in the solution:
\citet{RoscaniVoller2024} establish that for $\alpha<1$ the flux is piecewise
continuous \emph{with an infinite derivative at $s^-(t)$}, the divergence being
replaced at $\alpha=1$ by a finite derivative together with the classical jump.
A fixed grid must represent that divergence in the cell adjacent to the front, and
the undershoot is what that costs. The source itself is correspondingly singular:
it vanishes identically at $\alpha=1$ since $1/\Gamma(0)=0$, while the limit of
the $\alpha<1$ sources is a Dirac mass of unit weight. Deferring the source until a node lies a full cell behind the front,
carrying the running total, is \emph{worse} (about $500$ cells and triple the
front error) because the deferred total then arrives as one larger impulse.

\section{Verification}\label{sec:verify}

\subsection{Making the reference solution cheap}

\citet{KubicaRyszewska2021} give a self-similar solution with
$s(t)=c_1t^{\alpha/2}$, with $u$ a nested Neumann series built from
\begin{equation}
M_1(p,y)=\frac{1}{\Gamma(1-\alpha)}\int_y^p\!\left(1-p^{-2/\alpha}\mu^{2/\alpha}
\right)^{-\alpha}\mathrm{d}\mu .
\end{equation}
Substituting $w=(\mu/p)^{2/\alpha}$ reduces this to an incomplete Beta function,
\begin{equation}
M_1(p,y)=\frac{p\,\alpha}{2\Gamma(1-\alpha)}\,B\!\left(\tfrac{\alpha}{2},1-\alpha
\right)\left[1-I_{(y/p)^{2/\alpha}}\!\left(\tfrac{\alpha}{2},1-\alpha\right)
\right],
\end{equation}
agreeing with direct quadrature to $10^{-11}$ and making the endpoint singularity
analytic. The Neumann series is a Volterra composition whose discrete operator is
triangular and independent of $p$, so every row follows from one triangular solve
rather than a dense solve per row: $16.8$\,s to $0.046$\,s at $n=800$. Grid
convergence of the front constant is observed at order $1.39$ at $\alpha=0.6$,
consistent with the $2-\alpha$ rate implied by the $(1-w)^{-\alpha}$ singularity.

\subsection{Convergence of the scheme}

\begin{table}[t]
\caption{Front constant against the self-similar reference, $n_x=201$,
$n_t=600$, $\gamma=1$.}
\label{tab:verify}
\begin{tabular*}{\linewidth}{@{}llll@{}}
\toprule
$\alpha$ & exact $c_1$ & scheme & rel.\ error \\
\midrule
$0.5$ & $1.37514$ & $1.37421$ & $0.07\,\%$ \\
$0.7$ & $1.33929$ & $1.33958$ & $0.02\,\%$ \\
$0.9$ & $1.27889$ & $1.28033$ & $0.11\,\%$ \\
\bottomrule
\end{tabular*}
\end{table}

The exact solution has $s(0)=0$, but a fixed grid needs a finite
starting front $s_0$, which carries an $O(s_0^2)$ \emph{modelling} bias distinct
from discretisation error. Refining $\Delta x$ and $\Delta t$ at fixed $s_0$ drives
the spatial error below that bias and the total plateaus, resembling a
non-convergent scheme. Refining $s_0$ jointly restores clean second-order
behaviour: $5.5\times10^{-3}$, $1.5\times10^{-3}$, $4.4\times10^{-4}$,
$1.1\times10^{-4}$. The kick-start is a convergence parameter, and a refinement
study holding it fixed reports a spurious error floor.

\subsection{A traveling wave, and verification against a front that is not
self-similar}
\label{sec:tw}

Table~\ref{tab:verify} verifies the scheme against a self-similar solution
because that is the only kind the literature supplies. That is a real
restriction: the design's central claim is that it represents fronts no
front-fixing scheme can, and a self-similar benchmark cannot test it. An exact
solution with a genuinely non-self-similar front is available. To our knowledge
every previously published exact solution of this problem is self-similar, with
$s(t)\propto t^{\alpha/2}$~\citep{LiuXu2004,LiuXu2009,RoscaniBollatiTarzia2018,RoscaniTarzia2018,Ceretani2020,KubicaRyszewska2021,VollerRoscani2023,RoscaniVoller2024}.
\citet{Voller2010} solves a more general two-fractional-order problem, with
independent time order $\beta$ and space order $\alpha$, and finds the same
structure in a different exponent, $s(t)\propto t^{\beta/(\alpha+1)}$: a power
law, hence still self-similar. So are the solutions collected in the survey
of~\citet{RogosinDubatovskaya2023}, where every reported exact solution is expressed in Wright functions and the
conclusions name the construction of exact solutions as an open problem for
fractional Stefan-type problems; the word ``traveling'' does not occur in it. The
solution below is, to our knowledge, the first exact solution with the standard
Caputo derivative that is not self-similar, and the first traveling wave.
\citet{Garra2020} do obtain fronts that are not powers of $t$, in the
source-dominated limit, but only by replacing the Caputo derivative with a
derivative taken with respect to a chosen function of time $f(t)$, which gives
$s\propto f(t)^{\alpha/2}$; with the standard derivative, $f(t)=t$, their front
is again $t^{\alpha/2}$. Their memory also starts at $t=0$ rather than at
$h(x)$.

Seek a front of constant speed together with a profile depending on $x$ and $t$
only through the elapsed time since melting,
\begin{equation}
s(t)=s_0+vt,\qquad h(x)=\frac{x-s_0}{v},\qquad
u(x,t)=W(\tau)-1,\quad \tau=t-h(x)=\frac{s(t)-x}{v}.
\label{eq:tw}
\end{equation}
Then the source argument $t-h(x)$ is $\tau$ as well, $u_{xx}=W''/v^{2}$, and the
Caputo derivative with lower limit $h(x)$ becomes an ordinary one in $\tau$ with
lower limit $0$: \eqref{eq:consistent} collapses to a scalar fractional ODE.
Using $\Dc^{\alpha}V=\Drl^{\alpha}V$ for $V(0)=0$, and recognising the singular
source as exactly $-\Drl^{\alpha}[1]$, the whole problem is
\begin{equation}
\Drl^{\alpha}\!\left[W\right]=\frac{W''}{v^{2}},\qquad W(0)=1,\quad W'(0)=0 .
\label{eq:twode}
\end{equation}
The \emph{enthalpy} $W=u+1$ obeys a source-free equation: the singular source
disappears in this representation, as it does in the limit problem of
Section~\ref{sec:limit}.

Equation~\eqref{eq:twode} is a classical object. Differentiating it once more and
using $W'(0)=0$ puts it in single-term form, $\Dc^{\,2-\alpha}W=v^{2}W$ with
$2-\alpha\in(1,2)$: this is \citeauthor{Mainardi1996}'s fractional
\emph{oscillation} equation~\citep{Mainardi1996}, and its solution is the
one-parameter Mittag-Leffler function
\begin{equation}
W(\tau)=E_{2-\alpha}\!\left(v^{2}\tau^{\,2-\alpha}\right)
=\sum_{n\ge0}\frac{v^{2n}\,\tau^{\,n(2-\alpha)}}{\Gamma\!\left(n(2-\alpha)+1\right)} .
\label{eq:twseries}
\end{equation}
Matching a series in \eqref{eq:twode} directly reproduces
\eqref{eq:twseries}: the exponent lattice is forced, since two $x$-derivatives
must match one $\alpha$-derivative, and the resulting recursion telescopes,
because $\eta(\eta-1)\Gamma(n(2-\alpha)+1-\alpha)=\Gamma((n+1)(2-\alpha)+1)$
for $\eta=(n+1)(2-\alpha)$. Every property we need is then standard:
$E_\beta$ is entire, so \eqref{eq:twseries} converges for every $\tau$ and every
$v$; and its coefficients are positive, so $u=W-1\ge0$ automatically, in
agreement with the positivity theorem of~\citet{KubicaRyszewska2021}.

We stress that the profile equation is not the contribution. What is new is the
\emph{embedding}: that a traveling wave solves the free-boundary problem at all.
A time-fractional operator with the fixed lower limit $0$ is not
translation-invariant in $t$, so a traveling wave cannot satisfy it exactly: in
the time-fractional Fisher--KPP setting~\citet{Ishii2025} notes precisely this and
obtains fronts only asymptotically. The conservation-consistent formulation
escapes the objection because its lower limit $h(x)=s^{-1}(x)$ travels with the
front: in the moving frame the problem becomes autonomous, and the wave is exact.
Conservation consistency, adopted here on thermodynamic grounds, is what buys the
solution. In the classical limit the construction degenerates to the
constant-speed Stefan fronts of~\citet{Ricci1990}.

What makes \eqref{eq:tw} a solution of the free-boundary problem, rather than
merely of \eqref{eq:twode}, is that the two conditions not used to construct it
hold identically. First $u_x^-(s(t),t)=-W'(0)/v=0$, the emergent condition.
Second, since $W'(\tau)\sim\left[v^{2}/\Gamma(2-\alpha)\right]\tau^{1-\alpha}$
near the front while $\Drl^{1-\alpha}[\tau^{1-\alpha}]=\Gamma(2-\alpha)$, the
interface condition \eqref{eq:stefan} returns $\dot s=v$ exactly, \emph{for
every} $v$. The wave speed is free, and the boundary datum adjusts to sustain
it:
\begin{equation}
\gamma(t)=u(0,t)=W\!\left(t+s_0/v\right)-1 ,
\label{eq:twdatum}
\end{equation}
which varies in time and is not a power law. The benchmark therefore cannot be
run at all by a solver that admits only a constant datum.

Because \eqref{eq:twseries} was obtained by matching inside the bulk equation,
checking that equation again would be circular. We verify against the original
problem by three routes that do not use the recursion: the bulk residual, with
the Caputo term by adaptive quadrature from the definition and $u_{xx}$ by
fourth-order differences, is $\sim10^{-10}$ relative; the interface condition,
by numerical Riemann--Liouville differentiation, returns $v$ to between
$2\times10^{-9}$ and $3\times10^{-6}$ as $y\to s(t)^-$; and the global identity
$F(t)=F(0)-I^{\alpha}[u_x(0,\cdot)]$, which plays no part in the derivation,
holds to $1.3\times10^{-15}$. All three were run over
$\alpha\in\{0.3,0.5,0.7,0.85,0.9\}$ and $v\in\{0.5,1,1.5,2\}$.

The third check at first failed by $26\,\%$, because the identity assumes the boundary flux memory
begins at $t=0$, whereas a wave that has been travelling forever passed $x=0$ at
$h(0)=-s_0/v<0$. Restoring the correct lower limit took the defect from
$1.8\times10^{-1}$ to $5.6\times10^{-16}$, and identified the benchmark's one
structural cost. A fixed-grid solver started at $s_0$ under the convention
``lower limit $0$ and no source on $x<s(0)$'' discards the memory accumulated
over $[h(0),0]$. That discarded weight is $O(s_0^{2-\alpha})$, predicted from
$W'(\tau)\sim\tau^{1-\alpha}$ and measured at $1.692$, $1.494$, $1.298$, $1.104$
against $1.70$, $1.50$, $1.30$, $1.10$ for $\alpha=0.3,0.5,0.7,0.9$.

As with the self-similar benchmark, refining $\Delta t$ at fixed $s_0$
\emph{does not converge}. At $\alpha=0.5$ the
front error moves from $6.44\times10^{-3}$ to $6.33\times10^{-3}$ over a
fourfold refinement, $1.6\,\%$. The plateau is the pre-history floor, which
$\Delta t$ cannot touch. Refining $s_0$ jointly with the grid removes it, with
the measured slopes of Table~\ref{tab:tw}. We note that these exceed the
$2-\alpha$ figure above in every case and depend on $\alpha$ far more weakly, so
that estimate (derived from the discarded \emph{boundary} memory alone) is
conservative, plausibly because the source omitted on $x<s_0$ partly cancels it.
We have not proved a sharp rate, and the near-independence of $\alpha$ is
unexplained.

\begin{table}[t]
\caption{Traveling-wave benchmark, $v=1$, $T=0.3$. Left: convergence slope in
$s_0$ of the front error under joint $s_0$/grid refinement (Aitken-extrapolated
in $n_t$; raw $n_t=1600$ slopes in parentheses). Right: minimum field value with
\emph{no} grading active, against the bound of Section~\ref{sec:fix}.}
\label{tab:tw}
\begin{tabular*}{\linewidth}{@{}llll@{}}
\toprule
$\alpha$ & slope in $s_0$ & $2-\alpha$ & $\min u$ \\
\midrule
$0.3$ & $1.84$ $(1.83)$ & $1.70$ & $-7.3\times10^{-7}$ \\
$0.5$ & $1.76$ $(1.74)$ & $1.50$ & $-1.3\times10^{-5}$ \\
$0.7$ & $1.67$ $(1.65)$ & $1.30$ & $-8.4\times10^{-5}$ \\
$0.9$ & $1.61$ $(1.54)$ & $1.10$ & $-4.2\times10^{-4}$ \\
\bottomrule
\end{tabular*}
\end{table}

The right-hand column of Table~\ref{tab:tw} supports the diagnosis of
Section~\ref{sec:fix}. A constant-speed front makes the
front-locked time grid \emph{identical} to the uniform one, since
$\Delta t=\Delta x/v$ advances the front exactly one cell per step, and there is
no $t^{\alpha}$ start-up singularity for a uniform grid to stride past. If the
mechanism identified there is the right one, the positivity artifact should
nearly vanish here with all grading switched off, and it does, two to five
orders of magnitude below the $-0.10$ bound, against uniform-grid self-similar
values reaching $-0.98$. The final field is non-negative at every $\alpha$ tested.

\section{Measuring the inequivalence}\label{sec:noneq}

Normalising so the two formulations coincide identically at $\alpha=1$, their
front constants differ by $0.237$ at $\alpha=0.6$ ($17.4\,\%$ of the consistent
value) against a self-refinement discretisation bound of $6.3\times10^{-4}$, a
separation of nearly three orders of magnitude, closing monotonically as
$\alpha\to1$. This is~\citet{Ceretani2020}'s theorem made quantitative, and it is
why a scheme for one formulation cannot be validated against a similarity
solution of the other.

This figure needs a caveat. The consistent solver reproduces its own closed
form to $0.02\,\%$ at this $\alpha$, so the shared machinery is sound. The naive
solver does not: it returns $1.124$ against the closed form's $0.998$, a gap of
$12.6\,\%$, and under joint $s_0$/grid refinement it converges to the former.
The $17.4\,\%$ above is therefore the gap between the two \emph{solvers};
comparing the two closed forms instead gives $0.363$, or $26.6\,\%$. We have not
settled which reference is correct for the naive formulation. Testing its bulk
equation needs a value for $u(x,r)$ at times before the front reaches $x$,
which is the fixed-lower-limit ambiguity this paper is about. The
self-refinement bound measures grid convergence only, so it cannot decide
between the two figures. Both, however, are more than two orders of magnitude
above the discretisation error, which is the point: the formulations are not
equivalent below $\alpha=1$.

\section{The classical limit: a benchmark against a neural network}
\label{sec:pinn}

At $\alpha=1$ the obstruction of Section~\ref{sec:nocr} disappears, front fixing
becomes legitimate, and \eqref{eq:front} reduces to the classical balance. The
resulting solver is an ordinary implicit scheme of about forty lines. We use it
here for two purposes: as an independent check of the $\alpha\to1$ limit, and as a
comparator for a recent moving-boundary PINN.

Physics-informed neural networks minimise a PDE residual sampled inside the
domain together with a data or boundary loss~\citep{Raissi2019}; the fractional
case, where the residual itself requires evaluating a nonlocal operator, has its
own literature~\citep{Pang2019}. For moving boundaries specifically,
\citet{WangPerdikaris2021} train a network to represent the free boundary of a
classical Stefan problem directly. \citet{Malla2026} propose a two-network architecture (one network for the free
boundary, one for the field, with collocation points sampled inside the predicted
domain) and apply it to seven problem classes. That generality is their
principal contribution and is not in question here. Their Table~6 additionally
reports that the PINN is roughly ten times more accurate than an implicit
finite-difference method on the one-phase Stefan problem, at $380$\,s against
$2.3$\,s. An accuracy comparison of that kind is a statement about two
implementations, and for a problem with a closed-form solution and a front linear
in time a relative $L_2$ error of $1.5\times10^{-2}$ is not what finite
differences readily deliver.

We reimplemented their architecture as specified (eight hidden layers of widths
$(20,50,80,80,80,50,20,10)$, outputs formed as sums over the last hidden layer,
dynamic sampling, Adam with their schedule) and ran both methods on their
benchmark, for which $\theta=e^{t/10+1-x}$, $\zeta=1+t/10$ exactly.

\begin{table}[t]
\caption{The benchmark of~\citet{Malla2026}. Relative $L_2$ errors against the
exact solution. Rows marked ``reported'' are from their Table~6. Bold marks the
coarsest finite-difference configuration, the one compared against the PINN
median in the text.}
\label{tab:pinn}
\begin{tabular*}{\linewidth}{@{}llll@{}}
\toprule
Method & $L_2$ $\zeta$ & $L_2$ $\theta$ & time (s) \\
\midrule
PINN, reported             & $7.07\times10^{-4}$ & $1.42\times10^{-3}$ & $380.4$ \\
\quad ours, CPU, one seed  & $1.83\times10^{-4}$ & $5.71\times10^{-4}$ & $138.7$ \\
\quad ours, A30 GPU, median of 5 & $2.71\times10^{-4}$ & $6.27\times10^{-4}$ & $176.8$ \\
\quad \quad best / worst seed & $6.2\times10^{-5}$ / $1.06\times10^{-3}$ & --- & --- \\
FDM, reported              & $9.85\times10^{-3}$ & $1.51\times10^{-2}$ & $2.31$ \\
\quad ours, $n_x=51$, same host & $\mathbf{1.67\times10^{-4}}$ & $\mathbf{4.90\times10^{-4}}$ & $\mathbf{0.012}$ \\
\quad ours, $n_x=401$      & $2.59\times10^{-5}$ & $8.01\times10^{-5}$ & $0.646$ \\
\bottomrule
\end{tabular*}
\end{table}

Our reimplementation of their network reproduces and exceeds their reported
accuracy (Table~\ref{tab:pinn}), so the comparison is faithful and, if
anything, generous to it.

Two features of the replication matter for how such comparisons are
reported. First, training is not faster on a GPU: the median over five seeds on an NVIDIA A30 is $176.8$\,s against $138.7$\,s
on a single laptop CPU core. The networks and batches here are small enough that
kernel-launch overhead dominates, so the reported training cost is not an artefact
of CPU-bound execution. Second, the spread between seeds is a factor of seventeen, from
$6.2\times10^{-5}$ to $1.06\times10^{-3}$ in the front-position error.
A single training run is a draw from a wide distribution; our own first
single-seed result, $1.83\times10^{-4}$, ranks second of five and was a
favourable draw. We therefore report the median and range, and recommend
that others do the same.

Against the median, the finite-difference scheme at its coarsest useful resolution
is $1.6\times$ more accurate in $12$\,ms against $177$\,s on the same host (a
factor of $1.5\times10^{4}$ in time), and is $59\times$ more accurate than the
classical result reported alongside the network. Expressed as accuracy per unit
compute the gap is four orders of magnitude. The reported advantage of the network comes from a weak finite-difference
baseline, not from the network, consistent with~\citet{Grossmann2024}, who found finite elements ahead by two to six orders of
magnitude in one dimension on fixed domains.

\subsection{The obstruction, restated for neural solvers}

Moving-boundary PINNs frequently normalise the spatial coordinate, either
explicitly or by feeding a rescaled input to the network. Section~\ref{sec:nocr}
applies verbatim: for a fractional governing equation, residuals evaluated at
fixed $\xi=x/s(t)$ belong to a different problem, and the discrepancy is a
nonlocal integral rather than a correction term. The two agree at $\alpha\to1$,
which is precisely where a sanity check would be run and would therefore fail to
detect it. For fractional moving-boundary problems, residuals should be evaluated
on a fixed coordinate with the domain history represented explicitly (the same
constraint that shapes Section~\ref{sec:scheme}), or the correction must be
derived and included.

\section{The source-dominated limit, and where the difficulty really lies}
\label{sec:limit}

\citet{RoscaniVoller2024} obtain, as the vanishing-specific-heat limit of
\eqref{eq:consistent}, a problem in which the singular source is the \emph{only}
term balancing diffusion:
\begin{equation}
\frac{1}{\Gamma(1-\alpha)\,(t-h(x))^{\alpha}} = u_{xx},
\qquad s(0)=0,\ u(0,t)=1,\ u(s,t)=0,\ u_x(s,t)=0 .
\label{eq:limit}
\end{equation}
There is no time derivative: at each instant this is a two-point boundary value
problem carrying three conditions, the extra one determining $s(t)$. It is
therefore an unusually sharp test of a source treatment, since nothing else can
absorb an error, and it is precisely the regime in which our artifact lives.

We re-derive the exact solution rather than transcribe it. With
$\xi = x/(a t^{\alpha/2})$ and $h(x)=a^{-2/\alpha}x^{2/\alpha}$ one has
$t-h(x)=t(1-\xi^{2/\alpha})$, so \eqref{eq:limit} reduces to
$f''(\xi)=a^{2}[\Gamma(1-\alpha)(1-\xi^{2/\alpha})^{\alpha}]^{-1}$ with $f(0)=1$,
$f(1)=f'(1)=0$. Integrating twice from $\xi=1$ and exchanging the order of
integration,
\begin{equation}
f(\xi)=\int_\xi^1 (w-\xi)\,g(w)\,\mathrm{d}w,\qquad
g(w)=\frac{a^{2}}{\Gamma(1-\alpha)}\left(1-w^{2/\alpha}\right)^{-\alpha},
\end{equation}
and $f(0)=1$ forces $a^{2}\Gamma(\alpha+1)=2$, recovering their front constant.

\subsection{The same problem in two representations}

Marched in position space (finding at each step the $s$ for which
$\int_0^s w\,\sigma(w)\,\mathrm{d}w=1$, with $\sigma$ built from the recorded crossing
times), the scheme reproduces $a$ to $0.07\,\%$ at $\alpha\le0.5$ and then fails:
$5.8\,\%$ at $\alpha=0.7$ with a uniform rule, $40\,\%$ at $\alpha=0.9$. Grading
the quadrature towards the front rescues the moderate cases but destroys the
large-$\alpha$ ones, because clustering strongly enough to resolve
$(s-w)^{-\alpha}$ starves the smooth remainder of the interval. Both failures are
the same endpoint singularity that produces the undershoot of
Section~\ref{sec:source}.

Changing the integration variable from position to crossing time, $w=s(\tau)$ with
$h(w)=\tau$, removes it entirely:
\begin{equation}
\int_0^{s}\! w\,\sigma(w)\,\mathrm{d}w
=\frac{1}{2\Gamma(1-\alpha)}\int_0^t (s^{2})'(\tau)\,(t-\tau)^{-\alpha}\mathrm{d}\tau
=\tfrac{1}{2}\,\Dc^{\alpha}\!\left[s^{2}\right](t),
\end{equation}
so \eqref{eq:limit} is \emph{exactly} the scalar fractional differential equation
\begin{equation}
\Dc^{\alpha}\!\left[s^{2}\right](t)=2,\qquad s(0)=0,
\label{eq:fode}
\end{equation}
which $s^{2}=a^{2}t^{\alpha}$ solves with $a^{2}\Gamma(\alpha+1)=2$ as before.

\begin{table}[t]
\caption{Front constant for the limit problem \eqref{eq:limit}: the same problem
solved in position space and through the reduction \eqref{eq:fode}.}
\label{tab:limit}
\begin{tabular*}{\linewidth}{@{}llll@{}}
\toprule
$\alpha$ & exact $a$ & position space & via \eqref{eq:fode} \\
\midrule
$0.30$ & $1.492811$ & $0.07\,\%$ & $0.0026\,\%$ \\
$0.50$ & $1.502251$ & $0.07\,\%$ & $0.0043\,\%$ \\
$0.70$ & $1.483609$ & $5.8\,\%$  & $0.0058\,\%$ \\
$0.90$ & $1.442050$ & $40\,\%$   & $0.0050\,\%$ \\
$0.95$ & $1.428658$ & $>100\,\%$ & $0.0032\,\%$ \\
\bottomrule
\end{tabular*}
\end{table}

Table~\ref{tab:limit} makes the point. In \eqref{eq:fode} the offending
singularity is nothing but the ordinary L1 kernel, for which product-integration
weights are exact on the local polynomial, and four significant figures are
available uniformly in $\alpha$. The difficulty was never intrinsic to the
problem; it was a property of the representation.

We do not claim this disposes of the artifact in the full solver, where the
source is indexed by $h(x)$ pointwise in space and no such global reduction is
available. What it establishes is the mechanism and the direction: the undershoot
is the discrete cost of representing, on a fixed spatial grid, a singularity that
is benign when expressed in time. A front-adapted or time-indexed treatment of
the source term is the natural remedy, and \eqref{eq:fode} shows how much is
recoverable when one is available.

\subsection{An exact, non-self-similar benchmark}
\label{sec:nonselfsimilar}

The reduction \eqref{eq:fode} has a second use: it gives an exact solution
whose front is not a power of $t$. Nothing in the change of variable used $\gamma=1$: the constant on the
right of the front condition \emph{is} the boundary datum. With a time-varying
datum $u(0,t)=\gamma(t)$ the limit problem is therefore exactly
\begin{equation}
\Dc^{\alpha}\!\left[s^{2}\right](t)=2\gamma(t),\qquad s(0)=0,
\label{eq:fode-gamma}
\end{equation}
which is \emph{linear} in $s^{2}$. A datum assembled from powers then integrates
term by term under the Caputo power rule: for
$\gamma(t)=\sum_j c_j t^{q_j}$,
\begin{equation}
s(t)^{2}=2\sum_j c_j\,
\frac{\Gamma(q_j+1)}{\Gamma(q_j+1+\alpha)}\;t^{\,q_j+\alpha}.
\label{eq:nonselfsimilar}
\end{equation}

Equation~\eqref{eq:nonselfsimilar} is exact and closed-form (not a
perturbation, since it holds for arbitrary coefficients), and whenever two or
more powers are present the front is not self-similar: the local
log--log slope $\mathrm{d}\ln s/\mathrm{d}\ln t$ sweeps continuously from
$(q_{\min}+\alpha)/2$ to $(q_{\max}+\alpha)/2$, so $s\neq Ct^{m}$ for any $m$.
For $\gamma(t)=1+\tfrac32 t$ at $\alpha=0.6$ the slope runs from $0.300$ to
$0.679$ across $t\in[10^{-3},10]$, against the constant-datum case which holds
$0.300$ to twelve digits. Every exact solution otherwise available for this
problem (\citet{RoscaniVoller2024}'s Eq.~(78), the Neumann series of~\citet{KubicaRyszewska2021}, and our own re-derivation above) puts a
self-similar ansatz in by hand. The non-power-law fronts of~\citet{Garra2020}
come from changing the operator; this one comes from the boundary datum, with
the standard Caputo derivative unchanged.

Because \eqref{eq:nonselfsimilar} is derived \emph{through} the reduction,
verifying it with the reduction would be circular. We instead go back to the
original statement: take $s$ from \eqref{eq:nonselfsimilar}, obtain
$h(w)=s^{-1}(w)$ by root-finding, and evaluate
$\int_0^{s(t)}w\,q(w)\,\mathrm{d}w$ by adaptive quadrature. It returns
$\gamma(t)$ to between $10^{-10}$ and $10^{-15}$ across $\alpha\in\{0.3,0.6,0.85\}$
and coefficients of order unity.

Both limit solvers reproduce it at first order: the same rate they achieve on
the self-similar problem. The sharper of the two is the position-space march,
which never uses the reduction: it carries the crossing times $h(w)$ explicitly
and is the limit-problem analogue of the fixed-grid scheme of
Section~\ref{sec:scheme}. Against \eqref{eq:nonselfsimilar} at $\alpha=0.5$ its
final-time error halves on each refinement, $4.31\times10^{-4}\to
2.13\times10^{-4}\to1.06\times10^{-4}$, with no degradation relative to the
self-similar case. The fixed-grid design's ability to represent a front that is
not a power of $t$ is thus verified against an exact solution rather than
against a refined one. Section~\ref{sec:tw} carries the same check through the
full solver, by a different route: no reduction of the kind used here exists
there, but a traveling wave supplies an exact non-self-similar solution
directly.

\section{Two candidate remedies, a mechanism, and a resolution}
\label{sec:fix}

\begin{figure}[t]
\centering
\includegraphics[width=\linewidth]{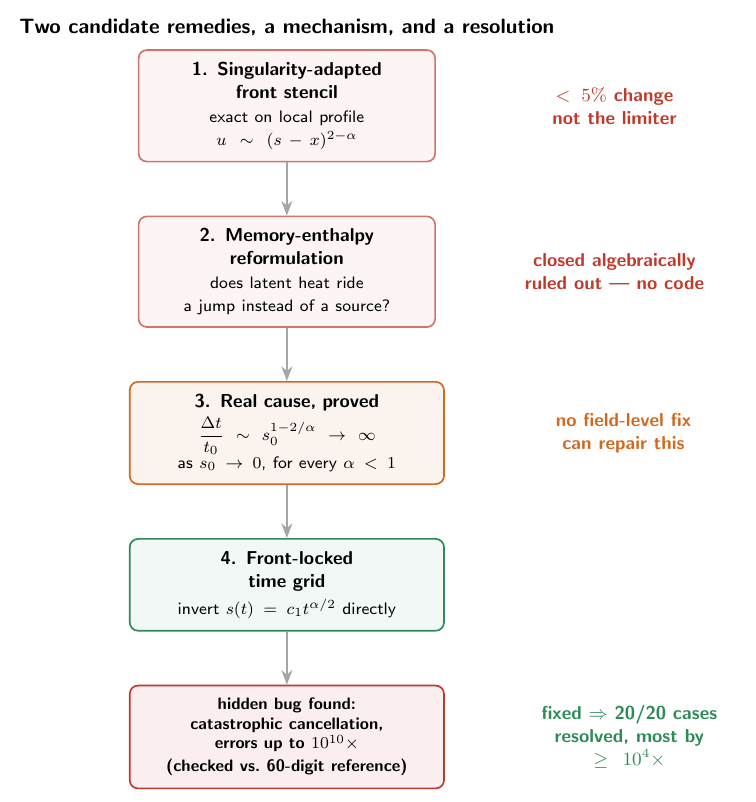}
\caption{The four approaches to the positivity violation, in the order tried.
The singularity-adapted stencil and the memory-enthalpy route are ruled out
(Section~\ref{sec:fix}); the scaling law~\eqref{eq:scaling} identifies the
cause (Section~\ref{sec:mechanism}); and the front-locked time grid, with the
non-uniform weights computed without cancellation, removes the violation in all
20 tested $(\alpha,s_0)$ cases (Section~\ref{sec:frontlocked}).}
\label{fig:remedyflow}
\end{figure}

We first tried two remedies for the positivity violation, and neither helped
(Fig.~\ref{fig:remedyflow}). A front-cell stencil fitted to the exact local
profile $u\sim(s-x)^{2-\alpha}$, on the basis $\{(s-x),(s-x)^{2-\alpha}\}$, is
exact where the quadratic stencil errs by $0.5\,\%$ at $\alpha=0.1$ rising to
$12\,\%$ at $\alpha=0.9$, yet under joint $s_0$/grid refinement it changes the
front-constant error by at most $5\,\%$ (ratios $1.00$, $1.00$, $1.04$, $1.05$
at $\alpha=0.3,0.5,0.7,0.9$). Switching to the memory-enthalpy formulation,
which~\citet{RoscaniVoller2024} prove equivalent, does not remove the source
either: their Eq.~(63) contains the same term,
$\Drl^{\alpha}_{h(x)}[cu+L]=c\,\Dc^{\alpha}_{h(x)}u+
L/[\Gamma(1-\alpha)(t-h(x))^{\alpha}]$, and the memory-enthalpy jump across the
front is zero (their Eqs.~46, 50). That route is closed without writing any
code.

\subsection{The mechanism}
\label{sec:mechanism}

The cause is the time grid. Under the joint $s_0$/grid refinement of
Table~\ref{tab:verify}, $\min u$ gets worse, not better: at $\alpha=0.5$ it
falls from $-6.3\times10^{-3}$ to $-6.0\times10^{-1}$ over the four rungs.
Refining the grid at fixed $s_0$ shrinks the violation, while shrinking $s_0$
on a fixed grid reproduces the blow-up on its own. In the worst case the front
moves $16.7$ times its own length in the first step and activates $133$ nodes
at once, far from the one-node-per-step regime that the source treatment of
Section~\ref{sec:source} assumes.

The exact front $s(t)=c_1t^{\alpha/2}$ reaches $s_0$ at
$t_0=(s_0/c_1)^{2/\alpha}$. Joint refinement takes $\Delta t\propto s_0$, so
\begin{equation}
\frac{\Delta t}{t_0}\ \sim\ s_0^{\,1-2/\alpha}\ \longrightarrow\ \infty
\qquad\text{as } s_0\to0,\ \text{for every } \alpha\in(0,1),
\label{eq:scaling}
\end{equation}
which grows like $s_0^{-3}$ at $\alpha=0.5$. The first time step jumps over the
whole interval in which the front reached $s_0$. Changing the initial field
does not help: seeding $u$ on $[0,s_0]$ with the exact profile at $t_0$ gives
$\min u=-6.021\times10^{-1}$, against $-5.982\times10^{-1}$ from a cold start
($\alpha=0.5$, $s_0=0.01$), because $t_0\sim10^{-9}$ while
$\Delta t\sim10^{-5}$. So on a uniform time grid, the front-constant
convergence of Table~\ref{tab:verify} and the positivity bound of
Table~\ref{tab:undershoot} hold in different parameter ranges.

\subsection{A front-locked time grid}
\label{sec:frontlocked}

Equation~\eqref{eq:scaling} points to the remedy: small time steps near $t=0$,
as in the graded meshes used for weakly singular time-fractional
problems~\citep{Stynes2017,Stynes2023}. A global power-law grid, $t_n=T(n/n_t)^r$, did not
work at any grading we tried, up to $r=2/\alpha$. What works is to build the
grid from the front law itself, so that the front moves exactly one cell per
step,
\begin{equation}
t_k=\Big(\tfrac{s_0+k\Delta x}{c_1}\Big)^{2/\alpha}-t_0,
\qquad k=0,1,2,\dots,
\label{eq:frontlocked}
\end{equation}
until the step reaches the uniform spacing, and uniform after that.

This grid needs the L1 Caputo weights and the Riemann--Liouville integral
of~\eqref{eq:front} on non-uniform steps. Both contain differences of nearly
equal powers, $(t_{n+1}-t_j)^{\kappa}-(t_{n+1}-t_{j+1})^{\kappa}$, whenever an
early, tiny step is weighed at a much later time. Computed directly, they lose
accuracy: against a 60-digit reference the Caputo weights are off by
$0.06$--$0.12\,\%$ each, which the history sum accumulates into a bias that
held the front-constant error near $4\,\%$, and the Riemann--Liouville weights
are off by up to ten orders of magnitude, sign included. We compute the Caputo
weights as $1-(1-r)^{\kappa}=-\operatorname{expm1}(\kappa\log 1p(-r))$ with
$r=(t_{j+1}-t_j)/(t_{n+1}-t_j)$, and the Riemann--Liouville weights with a
9-term Taylor series in $r$ below $r=10^{-2}$ and the closed form above it.
Both then match the 60-digit reference to machine precision.

\begin{table}[t]
\caption{Worst-case $\min u$ over $s_0\in\{0.005,\dots,0.08\}$, $n_x=801$,
$n_t=3200$: uniform-grid cold start against the front-locked grid
(Eq.~\eqref{eq:frontlocked}) with the cancellation fix, and the front-constant
error at that same $s_0$.}
\label{tab:graded}
\begin{tabular*}{\linewidth}{@{}llll@{}}
\toprule
$\alpha$ & uniform, cold & front-locked, fixed & front-const.\ error \\
\midrule
$0.3$ & $-0.452$ & $-1.9\times10^{-6}$ & $4.0\times10^{-4}$ \\
$0.5$ & $-0.979$ & $-4.0\times10^{-5}$ & $3.5\times10^{-4}$ \\
$0.7$ & $-0.463$ & $-1.6\times10^{-3}$ & $6.0\times10^{-5}$ \\
$0.9$ & $-0.035$ & $-1.1\times10^{-2}$ & $4.6\times10^{-4}$ \\
\bottomrule
\end{tabular*}
\end{table}

With these weights the front-locked grid removes the violation
(Table~\ref{tab:graded}). All twenty $(\alpha,s_0)$ cases clear the bound of
Table~\ref{tab:undershoot}, most by four orders of magnitude or more, and the
front-constant error stays at the level of the uniform-grid baseline, sometimes
slightly better and sometimes slightly worse. At $\alpha=0.3$, $s_0=0.01$,
previously the worst case, the front-constant error is $5.2\times10^{-4}$ and
$\min u=-1.4\times10^{-6}$.

The grid uses the exact front constant $c_1$, which the benchmark supplies and a
real problem does not. A practical version would estimate $c_1$ from a coarse
run, or build the grid as the front advances; we have tested neither. As
reported here, the front-locked grid is a diagnostic: it shows that the
positivity violation comes from the time grid, not from the spatial
discretisation or the source treatment.

\section{Conclusion}

Front fixing does not work for fractional moving-boundary problems, so we built
a scheme that does not need it. It uses a fixed grid, stores the time each node
melted, advances the front with an integrated energy balance, and times the
singular source to the sub-step crossing. Against an exact self-similar
solution the front constant is accurate to $0.02$--$0.11\,\%$. The competing
formulations differ by $17$--$27\,\%$ at $\alpha=0.6$, so a scheme must be
checked against a solution of its own formulation; the $\alpha\to1$ limit
cannot tell them apart.

We also gave two exact solutions whose fronts are not powers of $t$: a closed
form for the source-dominated limit with a time-varying boundary value, and a
Mittag-Leffler traveling wave for the full problem. To our knowledge these are
the first non-self-similar exact solutions of this problem with the standard
Caputo derivative, and the scheme
reproduces both.

The remaining positivity violation comes from the uniform time grid, through
the scaling law~\eqref{eq:scaling}. A grid locked to the front removes it once
the non-uniform weights are computed without cancellation.

In the classical limit the same solver, about forty lines long, is four orders
of magnitude ahead of a recent moving-boundary neural network in accuracy per
unit of compute. Neural solvers trained in a normalised coordinate meet the
same obstruction as front fixing when the derivative is fractional.

Two questions remain open. The front-locked grid needs the front constant in
advance, so for now it is a diagnostic rather than a general method. And the
gap between the formulations is $17.4\,\%$ or $26.6\,\%$ depending on which
reference is correct for the naive formulation, which turns on the same
fixed-lower-limit question this paper is about.

\section*{Reproducibility}
The scheme, the reference solution, the PINN reimplementation and every number
reported here are in the open \texttt{nlstefan} package, with tests asserting the
convergence orders and tolerances quoted.

\section*{Acknowledgements}
C.G.T.F.\ acknowledges the support of the Natural Sciences and Engineering
Research Council of Canada (NSERC) [RGPIN-2024-03989]. This research was enabled
in part by support provided by SHARCNET and the Digital Research Alliance of
Canada.

\bibliographystyle{plainnat}
\bibliography{refs}

\end{document}